%
\input amstex
\documentstyle{amsppt}
\def\C{{\Bbb C}}
\def\Z{{\Bbb Z}}

\def\O{{\Cal O}}

\def\sh#1#2{\hbox{${\Cal #1}  #2 $}}
\def\o#1{\operatorname{#1}}
\def\Hom{\o{Hom}}

\def\Ext{\o{Ext}}

\def\ch{\o{ch}}
\def\rk{\o{rk}}
\def\what#1{\widehat{#1}}
\def\iso{\kern.35em{\raise3pt\hbox{$\sim$}\kern-1.1em\to}
         \kern.3em}
\def\rest #1,#2{{#1}_{\vert #2}}
\def\proof{\demo{Proof}}
\catcode`\@=11
\def\cub{\vrule height 1.4ex width 1.4ex depth -.1ex}
\def\qed{\ifhmode\unskip\nobreak\fi\ifmmode\ifinner\else\hskip5
\p@\fi\fi
 \hfill\nobreak\hbox{\hskip5\p@\cub\hskip\p@}}
\def\whatx{\what X}

\def\hathsec#1{\widehat{\Cal #1}}
\def\hdhatp#1{R^{#1}\hat\pi_\ast}
\def\obligedskip{\vbox{\vskip4mm}}

\def\tr{\o{tr}}
\catcode`\@=\active
\magnification=1200
\pageno=1
\hyphenation{semi-stable}
\topmatter
\title  A FOURIER-MUKAI TRANSFORM FOR \\ STABLE BUNDLES
ON K3 SURFACES
\endtitle
\author  C\. Bartocci,$\star$
\ U\. Bruzzo,$\sharp$  and\
D\. Hern\'andez Ruip\'erez\P
\endauthor
\affil
$\star$\thinspace
Dipartimento di Matematica, Universit\`a di
Genova, Italia \\
$\sharp$\thinspace Scuola Internazionale Superiore di Studi
Avanzati \\ (SISSA --- ISAS), Trieste, Italia \\
\P\thinspace Departamento de Matem\'atica Pura y Aplicada, \\
Universidad de Salamanca, Espa\~na \\
\endaffil
\address $\star$
Dipartimento di Matematica, Universit\`a di
Genova, Via  Dodecaneso 35, 16146  Ge\-no\-va, Italy.
\endaddress
\email
 bartocci\@dima.unige.it
\endemail
\address
$\sharp$ S.I.S.S.A., Via Beirut 2-4, 34014 Miramare,
Trieste, Italy.
\endaddress
\email
bruzzo\@sissa.it
\endemail
\address
\P\ Departamento de Matem\'atica Pura y Aplicada,
Universidad de Salamanca, Plaza de la Merced 1-4, 37008
Salamanca, Spain.
\endaddress
\email
ruiperez\@gugu.usal.es
\endemail
\date
Revised --- 6 August 1996
\enddate
\thanks
Research partly supported  by the Italian Ministry for University and
Research through the research projects `Metodi geometrici
e probabilistici in fisica matematica' and `Geometria
delle  variet\`a differenziabili,'  by the Spanish
DGICYT through the research projects PB91-0188 and
PB92-0308,  by  EUROPROJ, and by the National
Group GNSAGA of the Italian Research
Council.
\endthanks
\abstract
We define a Fourier-Mukai transform for sheaves on K3 surfaces
over $\C$, and show that it maps polystable bundles to polystable
ones. The r\^ole of ``dual'' variety to the given K3 surface $X$
is here played by a suitable component $\whatx$ of the
moduli space of stable sheaves on $X$.
For a wide class of K3 surfaces
$\whatx$ can be chosen
to be isomorphic to $X$; then the Fourier-Mukai transform is
invertible, and the image of a zero-degree stable bundle $F$ is stable and
has the same Euler characteristic as $F$.
\endabstract
\endtopmatter
\leftheadtext{C\. Bartocci, U\. Bruzzo and D\. Hern\'andez Ruip\'erez}
\document
\head
1. Introduction and preliminaries
\endhead
Mukai's functor
may be defined within a fairly general setting; given two schemes $X$,
$Y$ (of finite type over an algebraically closed field $k$), and
an element $\Cal Q$ in the derived category $D(X\times Y)$ of $\O_{X\times
Y}$-modules,
Mukai \cite{18} defined a functor
${\Cal S}_{X\to Y}$ from the derived category
$D^-(X)$ to $D^-(Y)$,
$${\Cal S}_{X\to Y}({\Cal E}) = R\hat\pi_\ast
(\Cal Q{\overset L\to\otimes} \pi^\ast {\Cal E})$$ (here
$\pi: X\times Y\to X$ and $\hat\pi: X\times Y\to Y$ are the
natural projections).

Mukai has  proved that, when $X$ is an abelian variety,
$Y=\what X$ its dual variety, and $\Cal Q$ the Poincar\'e
bundle on $X\times\whatx$,   the functor  ${\Cal S}_{X\to \what X}$
gives rise to  an equivalence of categories.  If one is interested
in transforming sheaves rather than complexes one can
introduce (following  Mukai) the notion of WIT$_i$ sheaves:
an $\O_X$-module ${\Cal E}$ is said to be WIT$_i$ if its
Fourier-Mukai transform is  concentrated in degree
$i$, i.e\.  ${\Cal H}^kR\hat\pi_\ast (\Cal Q{\overset L\to\otimes}
\pi^\ast {\Cal E})=0$ for all $k\neq i$.

Let  $X$ is be a polarized abelian surface
and $\what X$ is its dual. If $\Cal F$ is a $\mu$-stable vector bundle  of
rank $\geq 2$ and zero degree, then it is
WIT$_1$, and its Mukai transform $\hathsec F$ is again a
$\mu$-stable bundle (with respect to a suitable polarization on $\whatx$) of
degree zero. This result  was proved  by Fahlaoui and Laszlo \cite{10}
and Maciocia
\cite{14}, albeit Schenk \cite{21} and  Braam and van Baal
\cite{7} had previously obtained a completely equivalent
result in the differential-geometric setting:  the Nahm-Fourier transform of
an instanton on a flat
four-torus (with no
flat factors) is an instanton over the dual torus (cf\. also
\cite{9},  and, for a detailed proof of the equivalence of the
two approaches, \cite{3}). The key remark which makes it
possible to relate the algebraic-geometric treatment to the
differential-geometric
one is that the Fourier-Mukai transform
is interpretable as the index of a suitable family of elliptic
operators parametrized by the target space $\what X$,  very much
in the spirit of Grothendieck-Illusie's approach to the
definition of  index of a relative elliptic complex.

The Fourier-Mukai transform can be studied also in the case of
nonabelian varieties. The general idea is to consider a variety
$X$, some moduli space $Y$ of  `geometric objects' over $X$, and
(if possible) a `universal sheaf'  $\Cal Q$ over the product
$X\times Y$. In the present paper we consider the case of a smooth
algebraic K3 surface $X$ over $\C$ with a fixed ample line
bundle $H$ on it; the r\^ole of dual variety is here played by the
moduli space $Y=M_H(v)$ of Gieseker-stable sheaves (with respect to $H$)
$\Cal E$ on $X$ having a fixed $v$, where
$$v({\Cal E})=\ch({\Cal E})\, \sqrt{\o{td}(X)}=(r,\alpha,s)\in H^0(X,\Z)\oplus
H^{1,1}(X,\Z)\oplus H^4(X,\Z)\,.$$
 This space, which is a  quasi-projective subscheme of
the moduli space of simple sheaves over $X$, has been extensively studied  by
Mukai \cite{19-20}. In particular, he proved the
following crucial result: if $M_H(v)= Y$ is nonempty and
compact of dimension two, then it is a K3 surface, isogenous to
the original K3 surface $X$.

If $X$ is a K3 surface satisfying these conditions and a few additional
assumptions, we can prove our first result (already announced in \cite{2}),
by exploiting the relationship between stable bundles and instanton
bundles (the so-called Hitchin-Kobayashi correspondence).

\proclaim{Theorem 1} Let $\Cal F$ be a
$\mu$-polystable IT$_1$ locally free sheaf of zero degree on $X$.
The Fourier-Mukai transform of $\Cal F$ is a $\mu$-polystable
locally free sheaf on $Y$.
\qed\endproclaim

For a wide family of K3 surfaces $X$,
we can choose
$v$ so that
$\whatx=M_H(v)$ is a K3 surface isomorphic to $X$.
The nonemptiness of $M_H(v)$ in this case has been proved in \cite{4} by
direct methods.
 Then we prove that the Fourier-Mukai transform is
invertible, just as in the case of abelian surfaces:
\proclaim{Theorem 2}
Let $\Cal F$ be an WIT$_i$ sheaf on $X$.
\sl Then its Fourier-Mukai transform
$\hathsec F=R^i\hat\pi_\ast(\pi^\ast\Cal F\otimes\Cal Q)$ is
a WIT$_{2-i}$ sheaf on $\whatx$,
whose Fourier-Mukai transform
$R^{2-i}\pi_\ast(\hat\pi^\ast\hathsec F\otimes\Cal Q^\ast)$ is
isomorphic to  $\Cal F$.
\qed\endproclaim
Furthermore, the Fourier-Mukai transform of a stable bundle is
stable:
\proclaim{Theorem 3}
Let $\Cal F$ be a zero-degree $\mu$-stable
bundle on $X$, with $v(\Cal F^\ast)\neq v$. Then $\Cal F$ is
IT$_1$, and its Fourier-Mukai transform
$\hathsec F$ is $\mu$-stable.
\qed\endproclaim

We end this section by fixing some terminology.
Let $X$ and $Y$ be compact complex manifolds, and let
$\Cal Q$ be a fixed coherent sheaf on  $X\times Y$, flat over $\O_Y$.
Let $\pi$, $\hat\pi$  be the projections onto the two factors of $X\times Y$.
\proclaim{Definition 1} A sheaf $\Cal E$ on $X$ satisfies the
$i$th Weak Index Theorem condition (i.e\. it is WIT$_i$)
if $R^j\hat\pi_\ast(\pi^\ast\Cal E\otimes\Cal Q)=0$ for
$j\neq i$; similarly, the sheaf $\Cal E$ is said to satisfy
the $i$th  Index Theorem condition (i.e\. it is IT$_i$) if
$H^j(\hat\pi^{-1}(y),\rest{\Cal E\otimes\Cal Q},{p^{-1}(y)})=0$
for $j\neq i$ and all $y\in Y$.
\endproclaim
The base change theorem implies that the sheaf
$\Cal E$ is IT$_i$ if and only if  it is both WIT$_i$ and
$R^i\hat\pi_\ast(\pi^\ast\Cal E\otimes\Cal Q)$ is locally free.

If either $\Cal Q$ is flat over $\O_{X\times Y}$, or $\Cal F$ is flat over
$\O_X$, the sheaves $R^k\hat\pi_\ast(\pi^\ast \Cal F\otimes\Cal Q)$ are the
cohomology sheaves of the Fourier-Mukai transform  $R\hat\pi_\ast(\pi^\ast
\Cal F{\overset L\to\otimes}\Cal Q)$  in the derived category. Then,
if $\Cal F$ is
WIT$_i$ we call the sheaf
$\hathsec F=R^i\hat\pi_\ast(\pi^\ast\Cal F\otimes\Cal Q)$ on
$Y$ its {\it Fourier-Mukai transform}.

\obligedskip
\head
2. Fourier-Mukai transform on K3 surfaces
\endhead

\subhead
2.1. Hyperk\"ahler manifolds and quaternionic instantons \endsubhead
A hyperk\"ahler manifold
is a $4n$-dimensional Riemannian
manifold $X$
which admits three complex
structures
$I$, $J$ and $K$, compatible with the Riemannian structure, such
that $IJ=K$.
On a hyperk\"ahler manifold one can introduce a generalized
notion of instanton \cite{15}.
The three
endomorphisms
$I$, $J$, $K$ of  $TX\otimes \C$ allow one to define an endomorphism
$\phi$ of
$\Lambda^2T^\ast X\otimes\C$,
$$\phi=I\otimes I+J\otimes J+K\otimes K\,.$$  This satisfies
$\phi^2=2\phi+3$,  so that at every $x\in X$
one has an
eigenspace decomposition
$$\Lambda^2(T_x^\ast X\otimes\C)=V_1\oplus V_2\tag 1
$$
corresponding to the eigenvalues 3 and $-1$ of $\phi$,
respectively. If $E$ is a $C^\infty$ complex vector bundle on $X$, with
connection
$\nabla$ and curvature
$R_\nabla$, we say that the pair  $(E,\nabla)$ is a
quaternionic instanton if $R_\nabla$, regarded as a section
of
$\o{End}(E)\otimes\Lambda^2T^\ast X$, has no component in $V_2$.
If $X$ has dimension 4 this  agrees with the usual
definition of instanton, since in that case the splitting (1) is
no more than the decomposition of the space of two-forms into
selfdual and anti-selfdual forms.
It is quite evident that $(E,\nabla)$ is a
quaternionic instanton if and only if the curvature $R_\nabla$
is of type (1,1) with respect to all the complex structures
of $X$ compatible with its hyperk\"ahler structure.
Moreover, $(E,\nabla)$ is an
Einstein-Hermite bundle with respect to all the induced K\"ahler
structures; thus, for any compatible complex structure on $X$, the bundle
$E$ admits a holomorphic structure, and the sheaf of its
holomorphic section is then $\mu$-polystable with respect to a polarization
given by the K\"ahler form determined by the given  complex structure
(we recall that a coherent sheaf is said to be $\mu$-polystable
if it is a direct sum of stable sheaves having the same slope).

\smallskip
\subhead
2.2. Assumptions on the moduli space
\endsubhead
 Let
$X$ be a projective K3 surface over $\C$. The cup product
defines a
$\Z$-valued pairing on the graded ring $H^{\bullet}(X,\Z)$:
$$(a,b,c)\cdot(a',b',c') =-ac' +bb'-ca'\,.$$
For every sheaf ${\Cal F}$ on $X$ we define the Mukai
vector $v({\Cal F})\in H^{\bullet}(X,\Z)$ as the element
$\operatorname{ch}(\Cal F)\sqrt{\operatorname{td}(X)}=(\rk{\Cal
F}, c_1({\Cal F}), s({\Cal F}))$, where
$s({\Cal F})=\rk{\Cal F}+\ch_2({\Cal F})=\chi(\Cal F)-\rk{\Cal F}$.

We fix a polarization $H$ on $X$ and denote by
$M_H(v)$ the moduli space of sheaves on $X$ which are
Gieseker-stable\footnote{We adopt the
usual definition of (semi)stability, so that a Gieseker- (or
$\mu$-) semistable sheaf is assumed to be torsion-free. Moreover, by
``(semi)stable'' we shall mean ``Gieseker-(semi)stable.''}
with respect to $H$, having a fixed
$v$.

We shall make the following assumption.

\smallskip\noindent{\bf A1.} \it  One can choose a Mukai vector
$v=(r,\ell,s)$ which is  primitive and isotropic and satisfies
 $\gcd$ $(r,\deg\ell,s)=1$. \rm
\smallskip
By results of
Maruyama \cite{16,17} and Mukai \cite{19,20} one has:
\proclaim{Proposition 1} If $v$ satisfies assumption A1,
then $M_H(v)$ is either empty or a projective
K3 surface and a fine moduli scheme parametrizing stable sheaves
with Mukai vector $v$, so that there is a universal sheaf $\Cal
Q$ on $X\times M_H(v)$.
\qed\endproclaim

Now we fix a Mukai vector $v=(r,\ell,s)$ satisfying the assumption A1
and  the following additional assumptions:
\smallskip\it
\noindent{\bf A2}. The divisor $\ell$ has degree zero and $r>1$;
\smallskip
\noindent{\bf A3}.
The moduli space $M_H(v)$ is not empty, and parametrizes $\mu$-stable
sheaves.\rm
\smallskip
In this situation the universal sheaf $\Cal Q$ on $X\times M_H(v)$
is locally free by virtue of
\cite{20}, Corollary 3.10, that we recall in the following form:

\proclaim{Proposition 2}  If $v=(r,\ell,s)$ is isotropic and $r>1$, every
$\mu$-stable sheaf ${\Cal F}$ on $X$ with $v({\Cal F})=v$ is
locally free.
\qed\endproclaim

\smallskip
\subhead
2.3. Universal sheaf and the Atiyah-Singer bundle
\endsubhead
We consider a K3 surface $X$ and a Mukai vector $v$ satisfying
assumptions A1, A2 and A3, so that the fine moduli space
$\whatx=M_H(v)$ (a K3 surface) is formed by locally free $\mu$-stable sheaves
of zero degree. According to
the Hitchin-Kobayashi correspondence,
these bundles correspond to irreducible $U(n)$-instantons, hence
$\whatx$ is identified with a moduli space $M$ of
instantons.

This identification can be described as follows.
The product manifold $X\times M$ carries
the universal Atiyah-Singer instanton bundle $P$ \cite{1},
which is a hermitian bundle with a universal hermitian connection.
The curvature of the latter is of type (1,1) with respect
to the natural complex structure of $X\times M$ \cite{12}.
Therefore, $P$ can be given a holomorphic structure, and
the resulting holomorphic vector bundle $\Cal P$ is relatively stable
with respect to the projection onto $M$, since the universal
connection is a family of instantonic connections. Then there is
a morphism $f\colon M\to\whatx$ such that
$(\operatorname{id}\times f)^\ast \Cal Q\otimes\hat\pi^\ast\Cal
L\simeq \Cal P$,
where $\Cal L$ is a holomorphic line bundle on
$M$, and $\hat\pi\colon X\times\whatx\to\whatx$ is the
projection. Since $\Cal Q$ is defined up to
tensoring by the pullback of a line bundle on $\whatx$
we may assume that $(\operatorname{id}\times f)^\ast\Cal Q\simeq \Cal P$.
The map $f$ is exactly the identification given by the
Hitchin-Kobayashi correspondence.

To any complex structure on $X$ there corresponds a well-determined
complex structure on $\whatx$. Then, if $X$ is hyperk\"ahler,
we may endow
the manifold $X\times\whatx$ with a natural
hyperk\"ahler structure.
\proclaim{Proposition 3}
The pair $(Q,\nabla)$, where $Q$ is the smooth bundle
underlying $\Cal Q$, and $\nabla$ is the universal connection,
is a quaternionic instanton.
\endproclaim
\proof
For any complex structure on $X$ compatible
with its hyperk\"ahler structure the curvature of the universal connection
is of type (1,1) \cite{12}.
\qed\enddemo

\subhead
2.4. Polystability of the Fourier-Mukai transform
\endsubhead
 We turn now to the main purpose of this paper, namely, the analysis
of the Fourier-Mukai functor for vector bundles on a K3 surface, and
prove that the Fourier-Mukai transform of a polystable bundle is
polystable.

Let $X$ and the Mukai vector $v$ satisfy the
assumptions A1,  A2 and A3, so that $X$ is a
projective K3 surface and $\whatx=M_H(v)$, a
moduli space of zero-degree  $\mu$-stable bundles on
$X$, is a  projective K3 surface.
Let $\pi$, $\hat\pi$ be the projections onto the two factors of
$X\times\whatx$.

We consider a $\mu$-polystable IT$_1$
locally free sheaf\footnote{As it will be
apparent in Section 4.4, the condition of being IT$_1$ is actually a
mild one.}
${\Cal F}$ of zero degree on $X$ and its Fourier-Mukai transform
$\hathsec F=\hdhatp1(\pi^\ast\Cal
F\otimes\Cal Q)$  on $\whatx$.
According to the Hitchin-Kobayashi correspondence,
$\Cal F$ is the sheaf of holomorphic  sections of a vector bundle
$F$ which carries an hermitian  metric, such that if $\nabla$ is the
connection determined by the holomorphic and hermitian
structures of $F$,  the pair $(F,\nabla)$ is an instanton.
The bundle $\pi^\ast F\otimes Q$ on $X\times\whatx$
may be endowed with a connection $\tilde\nabla$ obtained
by coupling the universal connection of $Q$ with the pullback connection
$\pi^\ast \nabla$, and the pair $(\pi^\ast F\otimes Q,\tilde\nabla)$
is a quaternionic instanton.

We can now state the main result of this section.
By regarding $\whatx$ as an
instanton moduli space, it carries a natural K\"ahler metric
$\Phi_{\whatx}$, usually called the {\it Weil-Petersson metric\/}.
\proclaim{Theorem 1} Let ${\Cal F}$ be a $\mu$-stable
IT$_1$ locally free sheaf of zero degree on $X$. The
Fourier-Mukai transform $\hathsec F$ is a
$\mu$-polystable locally free sheaf on $\whatx$ (with respect to
the K\"ahler form of $\whatx$).
\endproclaim
\demo{Proof}
Let us denote by $\frak C$ the set of complex structures
in $X$ which are compatible with its hyperk\"ahler structure, and shall write
$X_I$ for $X$ endowed with the complex structure $I$. For every choice of
$I$ the
moduli space $\whatx$ has an induced complex structure, so that $\frak C$
parametrizes also the complex structures in   $\whatx$ compatible with its
hyperk\"ahler structure. For every choice of $I$ we consider on the product
manifold the complex structure $X_I\otimes\whatx_I$, so that $X\times\whatx$
acquires a hyperk\"ahler structure, and  $\frak C$
parametrizes  the compatible complex structures.

Let $G$ denote the so-called isotropy group of $X$,
i.e\. the group which permutes the elements in $\frak C$. The action of
$G$  on $\frak C$ is transitive; we denote by  $g_{II'}$ an element in $G$
mapping
$I$ to $I'$.

For any choice of $I\in\frak C$ we have a relative Dolbeault
complex $\Omega^{0,\bullet}_I$ on $X\times \whatx\to\whatx$,
and a relative Dirac complex $\Cal S^\pm_I$, where
$$\Cal S^-_I=\Omega^{0,1}_I,\qquad S^+_I
=\Omega^{0,0}_I\oplus\Omega^{0,2}_I\,.$$
The relative Dirac operator $D_I$ acts as a $C^\infty_{\whatx}$-linear
morphism $\Cal S^+_I\to\Cal S^-_I$. The element $g_{II'}$ of the isotropy
group intertwines the Dirac operators associated with the two complex
structures, so that one has a commutative diagram
$$\CD
0 @>>> \Cal F^\infty_I @>>> \hat\pi_\ast(\widetilde{\Cal F}^\infty\otimes
\Cal S^-_I) @> D_I^\ast >>
\hat\pi_\ast(\widetilde{\Cal F}^\infty\otimes
\Cal S^+_I) @>>> 0 \\
@. @VVV @VV g_{II'}V    @VV g_{II'}V \\
0 @>>> \Cal F^\infty_{I'} @>>> \hat\pi_\ast(\widetilde{\Cal F}^\infty\otimes
\Cal S^-_{I'}) @> D_{I'}^\ast >>
\hat\pi_\ast(\widetilde{\Cal F}^\infty\otimes
\Cal S^+_{I}) @>>> 0
\endCD\ ;\tag 2  $$
here $D^\ast_I$ is actually the adjoint of the Dirac operator
associated with the complex structure $I$ coupled with the
connection of $\pi^\ast F\otimes Q$; moreover, $^\infty$ means that
we are considering the sheaf of smooth sections of a holomorphic bundle.
The operators $D^\ast$ are surjective due to the IT$_1$ condition, and
every  sheaf $\widehat{\Cal F}^\infty_I$ is the sheaf of smooth sections
of the Fourier-Mukai transforms $\widehat{\Cal F}_I$ (cf\. \cite{3,6}).

The diagram (2) induces an isomorphism $\hat g_{II'}\colon
\widehat F_I\to\widehat F_{I'}$
of $C^\infty$ vector bundles;
 moreover, since the bundle $\pi^\ast F\otimes Q$
has a natural hermitian metric, the horizontal arrows in this diagram allow
one to introduce an hermitian metric $\hat h_I$ on every bundle
 $\widehat F_{I}$, and $\hat g_{II'}$ is then an isometry.
By coupling
the relative connection induced by the connection on $\pi^\ast F\otimes Q$
with the relative Dolbeault operator associated with a complex structure $I$
and by taking direct images
one defines a connection $\hat\nabla_I$ on $\widehat F$. Since this
connection is also induced by the direct images of the coupled connections
on $\hat\pi_\ast(\widetilde{\Cal F}_I^\infty\otimes\Cal S^{\pm}_{I})$,
diagram (2) proves that $\hat g_{II'}$ transforms $\hat\nabla_I$ into
$\hat\nabla_{I'}$.
So actually one has
a single hermitian bundle $(\widehat F,\hat h)$ with a single connection
$\hat\nabla$ whose curvature
is of type (1,1) with respect to all compatible complex structures.

Then, the pair $(\widehat F,\hat \nabla)$ is an instanton.
As a consequence, the  Fourier-Mukai transform $\widehat{\Cal F}_I$
of $\Cal F_I$ is $\mu$-polystable with respect to the K\"ahler form over
$\whatx$ determined by $I$. \qed\enddemo

\obligedskip
\head
3. Nonemptiness of moduli spaces. Moduli spaces isomorphic to the surface.
\endhead
In this section we consider a wide class of K3 surfaces for which the
hypotheses  on $X$ and on the moduli space $\whatx$
stated in the previous section may be satisfied.
When $X$ is such a surface, $\whatx$
may be chosen so as to be isomorphic to $X$ itself, and
$X$ may be identified with a moduli space of stable bundles on $\whatx$,
in such a way that the relevant universal  sheaf is $\Cal Q^\ast$.
\smallskip
\subhead
3.1. Nonemptiness of the  moduli space
\endsubhead
Let $X$ be a K3 surface endowed with a
polarization
$H$ and a divisor $\ell$ such that $H^2=2$, $H\cdot\ell=0$ and
$\ell^2=-12$, so that $v^2=0$ with $v=(2,\ell,-3)$.
We shall call   the K3 surfaces satisfying
these assumptions {\it reflexive.}

Throughout the remaining part of this paper
we shall also make the further technical assumption that on $X$ there
are no nodal curves of degree 1 or 2 with respect to $H$.
This will hold generically. Indeed, the ample
divisor $H$ defines a double cover of ${\Bbb P}^2$ branched over
a sextic; the image of a nodal curve of degree 1 is a line
tritangent to the sextic, while
the image of a nodal curve of degree 2 is a conic, tangent
to the sextic at six points. Neither situation can arise in the
generic case \cite{5}.

\proclaim{Lemma 1} If $X$ is a reflexive K3 surface and $D\cdot
H>2$ for every nodal curve $D$, the divisor $E=\ell+2H$ is not
effective. Then, $H^i(X,\O_X(\ell+2H))=0$ for $i\ge 0$.
\endproclaim
\proof
Since $E^2=-4$, if $E$ is effective, it is not irreducible
and $E=D+F$ for some nodal curve $D$. Then $D\cdot H=3$ so that
$F\cdot H=1$ and $F$ is also irreducible. It follows that
$F^2\ge-2$. If $F^2\ge 0$, then $D\cdot F\le -1$, so that $D=F$
which is absurd. Thus, $F^2=-2$ and $F$ is a nodal curve of degree
1, a situation we are excluding.
\qed\enddemo

The following fundamental result is proved in \cite{4}.

\proclaim{Proposition 4} Let $X$ be a reflexive K3 surface such that
$\ell+2H$ is
not effective. The moduli space
$\whatx$ of stable sheaves
${\Cal E}$  with $v({\Cal E})=(2,\ell,-3)$
is not empty.  Moreover, every element in
$\whatx$ is $\mu$-stable, so that
a  reflexive K3 surface  meets assumptions A1, A2 and A3.
\qed\endproclaim
One should notice that the elements in $\whatx$   do not satisfy the lower
bounds on the discriminant established by Sorger\footnote{Talk given at the
Europroj Workshop ``Vector bundles and structure of moduli", Lambrecht
1994.} and Hirschowitz and Laszlo \cite{11}.

\smallskip
\subhead
3.2. The isomorphism $\whatx\iso X$
\endsubhead
We consider now a reflexive K3 surface $X$ satisfying the assumption on
nodal curves
described in the previous section, and show
that there is a natural isomorphism between
$X$ and $\whatx=M_H(v)$.

The following result is a direct consequence of Lemma 1.
\proclaim{Lemma 2}
 $\dim \Ext^1(\Cal I_p(\ell+2H),\O_X)=\dim H^1(X,\Cal I_p(\ell+2H))=1$ for
every point
$p\in X$, where ${\Cal I}_p$ is the ideal sheaf of $p$.
\qed\endproclaim

\proclaim{Lemma 3} Let $[\Cal E]\in
\whatx$. For every section of
$\Cal E(H)$ there is an exact sequence
$$ 0 @>>> \O_X @>>> \Cal E(H) @>>> {\Cal I}_p(\ell+2H) @>>> 0\,,
$$  where ${\Cal I}_p$ is the ideal sheaf of a point
$p\in X$. Moreover, $\dim H^0(X,\Cal E(H))=1$, so that the point $p$
depends only on
the sheaf $\Cal E$.
\endproclaim
\proof Since $H^2(X,\Cal E(H))=0$ and $\chi(\Cal E(H))=1$, the
sheaf $\Cal E(H)$ has at least one section, and
we have an exact sequence $0 @>>> \O_X @>>> \Cal E(H) @>>> {\Cal K} @>>> 0$.
By taking double duals we obtain the exact sequence
$$ 0 @>>> \O_X(D) @>>> \Cal E(H) @>>> {\Cal I}_Z\otimes{\Cal O}_X(\ell+2H-D)
@>>> 0\,,
$$
where $Z$ is a zero-dimensional subscheme and $D$ is an effective divisor
with $0\le
D\cdot H\le1$.

If $D\cdot H=1$, $D$ is an irreducible curve, so that $D^2\ge -2$.
Moreover,
$H-2D\not\equiv 0$, thus
Hodge index theorem implies that $0>(H-2D)^2=4D^2-2$ and we have
two cases, $D^2=0$ and $D^2=-2$.
If $D^2=0$, then $(H-2D)^2=-2$, so that either
$H-2D$ or $2D-H$ is effective, which is absurd. Thus, $D^2=-2$, and $D$ is
a nodal curve
of degree 1, a situation we are excluding.
Hence $D=0$ and $\o{length}(Z)=1$.

The second statement follows from the previous Lemma.
\qed\enddemo

These results imply that there exists a one-to-one map of sets
$\Psi\colon\whatx\hookrightarrow X$, given by $\Psi([\Cal E])=p$, where $p$ is
the point determined by Lemma 3. Our next aim is to prove that this map is
actually an isomorphism of schemes; to this end we need a result that
follows from
Grauert's base change theorem and Lemma 2.
\proclaim{Corollary 1} The sheaf $\O_X(H)$ is IT$_0$, and
its Fourier-Mukai transform
$\Cal N=\hat\pi_\ast(\Cal Q\otimes\pi^\ast\O_X(H))$ is a line
bundle on $\whatx$.
\qed\endproclaim
Now, the natural morphism
$\hat\pi^\ast\Cal N \to\Cal Q\otimes\pi^\ast\O_X(H)$ provides a section
$$ 0@>>>\O_{X\times\whatx}@>\sigma>>\Cal Q\otimes\pi^\ast\O_X(H)\otimes
\hat\pi^\ast\Cal N^{-1}@>>>\Cal K@>>>0\,.
$$
Let $j\colon Z\hookrightarrow X\times\whatx$ be
the closed subscheme of zeroes of $\sigma$, and let $p=\pi\circ
j\colon Z\to X$, $\hat p=\hat\pi\circ j\colon Z\to \whatx$ be
the proper morphisms induced by the projections $\pi$ and
$\hat\pi$.
\proclaim{Proposition 5} The morphism $\hat p\colon Z\iso\whatx$ is an
isomorphism of schemes and the map $\Psi$ is the composite morphism
$p\circ\hat
p^{-1}\colon\whatx\to X$. Moreover, $\Psi$ is an isomorphism of schemes.
\endproclaim
\proof One easily sees that for every (closed) point $\xi\in\whatx$,
$\sigma$ induces
a section $0@>>>\O_X@>\sigma_{\xi}>>\Cal E(H)@>>>\Cal K_{\xi}@>>>0$ of
$\Cal E(H)$.
By Lemma 3,
$\Cal K_{\xi}\simeq{\Cal
I}_{p(\xi)}(\ell+2H)$ for a well-defined point
$p(\xi)\in X$. Then, every closed fibre of
$\hat p\colon Z\to\whatx$  consists of a single point and $\hat p$ is a
proper finite
epimorphism of degree 1 by   Zariski's Main Theorem. Since
$\whatx$ is smooth,  $\hat p$ it is an isomorphism.
Moreover one has
$\Psi=p\circ\hat p^{-1}$. By Lemma 2 for every (closed) point
$\xi\in\whatx$ the fibre
$\Psi^{-1}(\Psi(\xi))$ is a single point. If $\Psi(\whatx)$ is the
scheme-theoretic
image of $\Psi$, $\what
X\to\Psi(\whatx)$ is a finite epimorphism of degree 1 as above, so that
$\dim\Psi(\whatx)=2$ and $\Psi(\whatx)=X$. The smoothness of $X$ yields
once more the result.\qed\enddemo

\proclaim{Corollary 2} Let $\Cal E$ be a sheaf which fits into an exact sequence
$$ 0 @>>> \O_X @>>> \Cal E(H) @>>> {\Cal I}_p(\ell+2H) @>>> 0\,,
$$  where ${\Cal I}_p$ is the ideal sheaf of a point
$p\in X$. Then $\Cal E$ is $\mu$-stable and locally free with $v(\Cal
E)=v=(2,\ell,-3)$ and
$\Psi([\Cal E])=p$.
\qed\endproclaim

\obligedskip
\head
4. Fourier-Mukai transform on reflexive K3 surfaces
\endhead

In this section we investigate the main properties of the Mukai
transform in the case of reflexive K3 surfaces satisfying the assumption on
nodal curves described in the previous section. For these K3 surfaces the
Fourier-Mukai functor is invertible, and Theorem 1
holds in a stronger form, in that the Fourier-Mukai transform
$\hathsec F$ of a stable bundle $\Cal F$ is itself stable.
We also prove the nice formula $\chi(\Cal F)=\chi(\hathsec F)$.
\smallskip \subhead
4.1. A natural polarization for $\whatx$
\endsubhead
If $[\Cal E]\in\whatx$, then $H^0(X,\Cal
E)=H^2(X,\Cal E)=0$ and $\dim H^1(X,\Cal E)=1$. It follows that $\O_X$ is
IT$_1$ and
that $\Cal O_X$ is a line bundle.
We can then normalize
$\Cal Q$  by twisting it by
$\hat\pi^\ast(\hdhatp1\Cal Q)^{-1}$, so that
$\hdhatp1\Cal Q\simeq\O_{\whatx}$.
We shall  henceforth assume that $\Cal Q$ is normalized in this way.

Let us denote $\gamma=\ch(\Cal Q)$. If $\gamma^{i,j}$ is the
$(i,j)$ K\"unneth component of $\gamma$, one has $\gamma^{2,0}=\ell$; we
set $\gamma^{0,2}=-\hat\ell$. Then, Riemann-Roch theorem gives
$(-1)^i c_1(\widehat{\Cal L})=-\ch_2(\Cal
L)\,\hat\ell+\hat\pi_\ast(\gamma^{2,2}\, c_1(\Cal L))$
for every WIT$_i$ line bundle $\Cal L$. In particular, since $\O_X(H)$ is IT$_0$
and
$\widehat{\O_X(H)}$ is a line bundle (Corollary 1), one has
$$
c_1(\widehat{\O_X(H)})=-\hat\ell-\what H\,,\qquad \widehat{\O_X(H)}\iso
\O_{\whatx}(-\hat\ell-\what H)\,.\tag 3
$$
where $\what H=-\hat\pi_\ast(\gamma^{2,2}\,H)$.

We will show that the divisor $\what H$ is a natural polarization on
$\whatx$. Indeed,
the space
$\whatx$, regarded as a moduli space of instantons on
$X$, carries the Weil-Petersson metric $\Phi_{\whatx}$, already
considered in Theorem 1, and we can prove, in the  spirit of \cite{13},
that the class of this metric may be identified with the class
$\what H$.

\proclaim{Proposition 6}
$\what H=\left\lbrack\frac1{8\pi^2}\Phi_{\whatx}\right\rbrack$.
\endproclaim
\proof
Let ${\frak R}$ denote the  curvature of
the universal connection on $Q$ (cf\. Proposition 3).
Let $\Phi_X$ be the K\"ahler form of $X$, so that $\lbrack\Phi_X
\rbrack=H$.
We note the identities
$$\int_X\Phi_X\wedge(\tr{{\frak R}})^2=
\left(\frac{2\pi}{i}\right)^2
\int_X\Phi_X\wedge(c_1(\Cal Q))^2=
2(4\pi)^2(\ell\cdot H)\hat\ell=0$$
$$\int_X\Phi_X\wedge\tr{\frak R}^{2,0}\wedge{\frak R}^{0,2}=0\,.$$
By representing the Chern character $\gamma$ in terms of the
curvature form $\widetilde{\frak R}$ we may compute
$$\align \what H &=\frac1{8\pi^2}\int_X\Phi_X\wedge\tr{\frak R}^2 \\
&=\frac1{8\pi^2}\int_X\Phi_X\wedge\tr(2{\frak R}^{2,0}\wedge{\frak R}^{0,2}+
{\frak R}^{1,1}\wedge{\frak R}^{1,1}) =\frac1{8\pi^2}\Phi_{\whatx}\,.\endalign$$
\qed\enddemo
Thus $\what H$ is ample, and can be taken as a polarization
on $\whatx$; moreover, $\hat\ell\cdot\what H=0$, and
$\hat\ell^2=-12$, so that $\hat v=(2,\hat\ell,-3)$ is an isotropic Mukai
vector and $\what X$ is  a reflexive K3 surface with respect to $(\what
H,\hat\ell)$.

\smallskip
\subhead
4.2. $X$ as a moduli space of bundles on $\whatx$
\endsubhead
In this section we prove that we can regard $\whatx$ as a moduli space of
bundles with topological invariants $(2,\hat\ell,-3)$  that
are $\mu$-stable with respect to the natural polarization $\what H$; it
turns out that the relevant universal bundle in this case is
simply $\Cal Q^\ast$.

Lemma 3 suggests that the universal sheaf $\Cal Q$ can be obtained as an
extension of suitable torsion-free rank-one sheaves on $X\times\whatx$.
Let ${\Cal I}_\Psi$ be the ideal sheaf
of the graph
$\Gamma_\Psi\colon \whatx\hookrightarrow X\times\whatx$ of
$\Psi$.
\proclaim{Lemma 4} The direct image
$\hat\pi_\ast[\sh Ext^1({\Cal I}_\Psi\otimes\pi^\ast\O_X
(\ell+2H),\O_{X\times\whatx})]$ is a line bundle $\Cal L$ on
$\whatx$.
\endproclaim
\proof
Write $E=\ell+2H$ and $\O_\Psi=(\Gamma_\Psi)_\ast\O_{\whatx}$. Then,
$$
\sh
Ext^1({\Cal I}_\Psi\otimes\pi^\ast\O_X (E),\O_{X\times\whatx})\simeq
\O_\Psi\otimes\pi^\ast\O_X(-E)\,.
$$
By Lemma 1, $R^i\hat\pi_\ast\pi^\ast\O_X(-E)=0$
for $i\ge 0$, hence, from the exact sequence
$$
0@>>>{\Cal I}_\Psi\otimes\pi^\ast\O_X (-E)@>>>\pi^\ast\O_X (-E)@>>>
\O_\Psi\otimes\pi^\ast\O_X(-E)@>>>0\,,
$$
we obtain  $\hat\pi_\ast(\O_\Psi\otimes\pi^\ast\O_X(-E))\iso
R^1\hat\pi_\ast({\Cal I}_\Psi\otimes\allowmathbreak
\pi^\ast\O_X(-E))$. But for every
$\xi\in\whatx$ one has
$H^1(X,{\Cal
I}_\Psi\otimes\pi^\ast\O_X(-E)\otimes\kappa(\xi))=H^1(X,\Cal I_p(-E))$, where
$p=\Psi(\xi)$, and one concludes by Lemma 1 and by Grauert's base-change
theorem.
\qed\enddemo
It follows that the sheaf $\sh Ext^1({\Cal I}_\Psi\otimes\pi^\ast\O_X
(\ell+2H),\hat\pi^\ast(\Cal L^{-1}))$ has a section, so that there is an
extension
$$
0@>>>\hat\pi^\ast(\Cal L^{-1})@>>>\Cal P@>>>{\Cal I}_\Psi\otimes\pi^\ast\O_X
(\ell+2H)@>>>0\,.\tag 4
$$
Moreover, Lemma 3 implies that $\Cal P\otimes\pi^\ast\O_X(-H)$ is a
universal sheaf
on $X\times\whatx$; thus $\Cal P\iso\Cal
Q\otimes\hat\pi^\ast\Cal N\otimes\pi^\ast\O_X(H)$ for a line bundle $\Cal
N$ on $\whatx$.
The sheaves $\Cal L$ and $\Cal N$ are readily determined; by applying
$\hat\pi_\ast$ to the  sequence above one obtains
$$
\Cal L^{-1}=\widehat{\O_X(H)}\otimes\Cal N\iso
\O_{\whatx}(-\hat\ell-\what H)\otimes\Cal N\,,
$$
where the second equality is due to equation (3). Now,  by restricting the exact
sequence (4) to a fibre $\pi^{-1}(p)$, we obtain
$c_1(\Cal N)=-\hat\ell-\what H-c_1(\rest{\Cal
Q},{\pi^{-1}(p)})=-\what H$. Then we have
\proclaim{Proposition 7}
The sequence of coherent sheaves on
$X\times\whatx$
$$ 0@>>>\hat\pi^\ast\O_{\whatx}(-\hat\ell-2\what H)@>>>\Cal
Q\otimes\hat\pi^\ast\O_{\whatx}(-\what
H)\otimes\pi^\ast\O_X(H)@>>>
{\Cal I}_\Psi\otimes\pi^\ast\O_X
(\ell+2H)@>>>0\tag 5
$$
is exact.
\qed\endproclaim

This sequence allows us to compute the Chern character of
$\Cal Q$. In particular, we obtain
$$
\gamma^{2,2}=(\ell+2H)\cup \what H +H\cup\hat\ell-\iota\,,\tag 6
$$
where $\iota\in H^2(X,\Z)\otimes  H^2(\whatx,\Z)$ is the element
corresponding to the isomorphism $\Psi^\ast\colon$ $H^2(X,\Z)$ $\to
H^2(\whatx,\Z)$. From this we get
$$ \what H = \Psi^\ast(2\ell+5H),\qquad
\hat\ell=\Psi^\ast(-5\ell-12 H)\,.$$ 
By taking duals in the sequence
(5) and restricting to the fibres of $\hat\pi$ we obtain
$$ 0@>>>\O_{\whatx} @>>>
\Cal Q^\ast_p(\what H) @>>> {\Cal
I}_\xi (\hat\ell+2\what H) @>>>0\,,\tag 7 $$
where $p=\Psi(\xi)$.

We need to show that the sheaves $\Cal Q^\ast_p$ are $\mu$-stable
with respect to $\what H$. We are not in a position to apply
Corollary 4 to the reflexive K3 surface $(\whatx,\what
H,\hat\ell)$, since we cannot a priori exclude that $\whatx$
contains nodal curves of degree 1   with respect to
$\what H$. This problem
is circumvented as follows.
Since $\hat\ell+2\what
H=-\Psi^\ast(\ell+2H)$, it has negative degree with respect to $\Psi^\ast
H$, so that
it is not effective. Thus,
$\dim H^0(X,\Cal Q^\ast_p(\what H))=1$ and
$\dim\Ext^1({\Cal I}_\xi (\hat\ell+2\what H),\O_{\whatx})=1$.
An easy calculation shows that the sheaves $\Cal Q^\ast_p$ are simple,
so that $\Cal Q^\ast$ defines a morphism  $$
\alpha\colon X\to \o{Spl}(\hat v,\whatx)
$$
into the moduli scheme of simple sheaves on $\whatx$ with Mukai vector
$\hat v$.
Proceeding as in the proof of Proposition 5, we obtain that $\alpha$
is an isomorphism with a connected component of $\o{Spl}(\hat v,\whatx)$.
Moreover, since $\hat\ell+2\what
H$ is not effective, Proposition 4 for $(\whatx,\what H,\hat\ell)$ implies
that the
moduli space $M_{\what H}(\hat v)$ of stable sheaves on $\whatx$ (with respect
to $\what H$) with Mukai vector
$\hat v$ is a non-empty connected component of
$\o{Spl}(\hat v,\whatx)$, consisting of locally free $\mu$-stable sheaves.

According to the proof of Lemma 3, if $[\Cal F]\in M_{\what H}(\hat
v)$ the sheaf $\Cal F$ fits into an exact sequence like (7) for a
well-defined point
$\xi\in\whatx$ unless
$\Cal F$ is given by an extension
$$
0@>>>\O_{\whatx}(D)@>>>\Cal F(\what H)@>>> \Cal I_Z(\hat\ell+2\what H-D)@>>>0
$$
where $D$ is a nodal curve with $D\cdot\what H=1$ and $Z$ is a zero-dimensional
closed subscheme of $\what X$. In the latter case,
$H^i(\whatx,\Cal I_Z(\hat\ell+2\what H-D))=0$ for $i\ge 0$ so that
$\o{lenght}(Z)=0$ and $-4=(\hat\ell+2\what H-D)^2$. Then, $D\cdot
\hat\ell=-3$ and $D\cdot \Psi^\ast H=D\cdot(5\what H+2\hat\ell)=-1$, which
is absurd since $\Psi^\ast H$ is ample. Then, one has
$$ 0@>>>\O_{\whatx} @>>>
\Cal F(\what H)@>>> {\Cal
I}_\xi (\hat\ell+2\what H) @>>>0\,, $$
for a point $\xi\in\whatx$, and $\Cal F\iso \Cal Q^\ast_p$ with
$p=\Psi(\xi)$, since
$\dim\Ext^1({\Cal I}_\xi (\hat\ell+2\what H),\O_{\whatx})=1$.
Thus,  $M_{\what H}(\hat v)$ is contained in $\alpha(X)$ and the two spaces
must coincide.
This means that the bundles
$\Cal Q^\ast_p$ are $\mu$-stable with respect to $\what H$; therefore the
sequence (5)
exhibits explicitly the parametrization of vector bundles on
$\whatx$ with invariants $(2,\hat\ell,-3)$ that are $\mu$-stable with
respect to $\what H$ by the points of $X$. As a consequence:
\proclaim{Proposition 8}
$X$ is a fine moduli space of $\mu$-stable bundles on $\whatx$
(polarized by $\what H$) with invariants $(2,\hat\ell,-3)$, and the
relevant universal sheaf is $\Cal Q^\ast$.
\qed\endproclaim

\smallskip
\subhead
4.3. Inversion of the Fourier-Mukai transform
\endsubhead
Let $X$ be a K3 surface and $v$ a Mukai vector satisfying
assumptions A1, A2 and A3.
We consider  the Fourier-Mukai transform as a functor
$$
\Cal S_X(F)=R\hat\pi_\ast(\pi^\ast F{\overset L\to\otimes}\Cal Q)
$$ between the derived categories $D(X)$ and $D(\whatx)$
(we may use the full derived categories instead of the $D^-$ categories
because $\Cal Q$ is locally free). In view of Proposition 8 a
natural candidate for the inverse of $\Cal S_X$ is the functor
$\Cal S_{\whatx}\colon D(\whatx)\to D(X)$
$$
\Cal S_{\whatx}(G)=R\pi_\ast(\hat\pi^\ast G{\overset L\to\otimes}\Cal
Q^\ast)\,.
$$

Since $X$ and $\whatx$ are K3 surfaces, the relative dualizing complexes
$\omega_\pi$
and $\omega_{\hat\pi}$ are both isomorphic to $\O_{X\times\whatx}[2]$.
Then, a straightforward application of relative duality gives:

\proclaim{Proposition 9}
For every objects $F$ in $D(X)$ and $G$ in
$D(\whatx)$,
we have functorial isomorphisms
$$
\align
\Hom_{D(\whatx)}(G,S_X(F))&\iso\Hom_{D(X)}(S_{\whatx}(G),F[-2])\\
\Hom_{D(X)}(F,S_{\whatx}(G))&
\iso\Hom_{D(\whatx)}(S_X(F),G[-2])\,.
\endalign
$$
\qed\endproclaim

\proclaim{Proposition 10}
For every $G\in D(\whatx)$ there is a functorial isomorphism
$$
\Cal S_X(\Cal S_{\whatx}(G))\iso G[-2]
$$ in the
derived category $D(\whatx)$. Moreover, if $X$ is reflexive and $D\cdot H>2$ for
every nodal curve $D$ in $X$, then for every $F\in D(X)$ there is also a
functorial isomorphism
$$
\Cal S_{\whatx}(\Cal S_X(F))\iso F[-2]
$$ in the derived category $D(X)$.
\endproclaim
\proof  Let $q_1$ and $q_2$ be the projections onto
the two factors of
$\whatx\times\whatx$, and $\pi_{ij}$ the projection of
$X\times\whatx\times\whatx$ onto the product of the $i$th and
$j$th factors. Then, the composite functor is given by
$\Cal S_X(\Cal S_{\whatx}(G))=R q_{2,\ast}(q_1^\ast
G{\overset L\to\otimes}\widetilde{\Cal Q})$ (see \cite{18}), with
$$
\widetilde{\Cal Q}=R\pi_{23,\ast}(\pi_{12}^\ast\Cal
Q^\ast{\overset L\to\otimes}\pi_{13}^\ast\Cal Q)\iso R\pi_{23,\ast} R\sh
Hom^{\sssize{\bullet}} (\pi_{12}^\ast\Cal Q,\pi_{13}^\ast\Cal
Q)\,,
$$
where $R\sh Hom^{\sssize{\bullet}}(\ ,\ )$ denotes the total
derived functor of the complex  of sheaf homomorphisms.
By \cite{20}, Proposition 4.10, the right-hand side is
isomorphic in the derived category to
$\delta_\ast\Cal M[-2]$, where
$\delta\colon\whatx\hookrightarrow \whatx\times\whatx$ is the
diagonal embedding and $\Cal M$ is an invertible sheaf on
$\whatx$. It follows that
$\Cal S_X(\Cal S_{\whatx}(G))\simeq G{\overset L\to\otimes}\Cal M[-2]$. Let
$\rho\colon\whatx\times\whatx\to\whatx\times\whatx$ be the permutation
morphism; by
base-change theory we have
$$
\rho^\ast
\widetilde{\Cal Q}\simeq R\pi_{23,\ast}(\pi_{13}^\ast\Cal
Q^\ast{\overset L\to\otimes}\pi_{12}^\ast\Cal Q)
\simeq R\pi_{23,\ast}[(\pi_{13}^\ast\Cal
Q{\overset L\to\otimes}\pi_{12}^\ast\Cal
Q^\ast)^\ast]\simeq\widetilde{\Cal Q}^\ast[-2]\,,
$$ by relative duality for $\pi_{23}$. From $\rho\circ\delta=\rho$ one finds
that $\Cal M[-2]\simeq\Cal M^\ast[-2]$. Then, there is an isomorphism of
invertible sheaves $\Cal M\simeq\Cal M^\ast$, and  $\Cal M\simeq\Cal
O_{\whatx}$.

The second statement follows from the first by Proposition 8.
\qed\enddemo

\proclaim{Corollary 3} If $X$ is reflexive and $D\cdot H>2$ for
every nodal curve $D$ in $X$, there are functorial isomorphisms
$$
\align
\Hom_{D(\whatx)}(G,\bar G)&\iso
\Hom_{D(X)}(S_{\whatx}(G),S_{\whatx}(\bar G))\\
\Hom_{D(X)}(F,\bar F)&\iso
\Hom_{D(\whatx)}({\Cal S}_X(F),{\Cal S}_X(\bar F))\,,
\endalign
$$
for $F$, $\bar F$ in $D(X)$ and $G$, $\bar G$ in
$D(\whatx)$.
\qed\endproclaim

Thus,
there is a duality between the varieties $X$ and $\whatx$, in a
complete analogy with the case of the Fourier-Mukai transform on
abelian surfaces.

In particular we have the following result.
\proclaim{Theorem 2} Let $\Cal F$ be a WIT$_i$ sheaf on
$X$.  Then its Fourier-Mukai transform
$\hathsec F=R^i\hat\pi_\ast(\pi^\ast\Cal F\otimes\Cal Q)$ is
a WIT$_{2-i}$ sheaf on $\whatx$,
whose Fourier-Mukai transform
$R^{2-i}\pi_\ast(\hat\pi^\ast\hathsec F\otimes\hathsec Q^\ast)$ is
isomorphic to $\Cal F$.
\qed\endproclaim
\smallskip
\subhead
4.4. Stability of the Fourier-Mukai transform
\endsubhead
 We may
now prove a stronger form of Theorem
1. Let $X$ be a reflexive K3 surface such that
 $D\cdot H>2$ for every nodal curve $D$,
and choose $\whatx$ as before.
By proceeding as in Corollary 2.5 of \cite{18} and taking into account
Corollary 3, we obtain:
\proclaim{Lemma  5} Let $\Cal F$, $\Cal F'$ be coherent
sheaves on $X$. If $\Cal F$ is WIT$_i$ and $\Cal F'$ is
WIT$_j$, we have
$$
\Ext^h(\Cal F,\Cal F')\simeq\Ext^{h+i-j}(\hathsec F,\hathsec F')\,.
$$
for $h=0,1,2$. In particular, there is an isomorphism
$\Ext^h(\Cal F,\Cal F)\simeq\Ext^h(\hathsec F,\hathsec F)$
for every $h$, so that $\hathsec F$ is simple for every simple
WIT$_i$ sheaf $\Cal F$.
\qed\endproclaim

\proclaim{Lemma 5} If $\Cal F$ is a $\mu$-stable vector bundle
of degree zero, and $v(\Cal F^\ast)\neq (2,\ell,-3)$,  then $\Cal F$ is
IT$_1$. \endproclaim
\proof For every $\xi\in\whatx$ we have
 $H^2(X,\Cal F\otimes\Cal Q_\xi)^\ast\simeq\Hom(\Cal
Q_\xi,\Cal F^\ast)$. Since $\Cal F$ and $\Cal Q_\xi$ are
$\mu$-stable,  if there exists a nonzero morphism $\Cal
Q_\xi\to\Cal F^\ast$, then it is an isomorphism, which is
incompatible with the condition in the statement. The same
argument also shows that $H^0(X,\Cal F\otimes\Cal
Q_\xi)\simeq\Hom(\Cal F^\ast,\Cal Q_\xi)=0$, thus concluding the
proof.
\qed\enddemo

\proclaim{Theorem 3} Let $\Cal F$ be a zero-degree
$\mu$-stable bundle on $X$, with $v(\Cal F^\ast)\neq
(2,\ell,-3)$. Then its Fourier-Mukai transform
$\hathsec F$ is $\mu$-stable.
\endproclaim
\proof $\hathsec F$ is $\mu$-polystable by Theorem
1, and  simple by  Lemma 5,
so that it is $\mu$-stable.
\qed\enddemo

\smallskip
\subhead
4.4. Topological invariants
\endsubhead

We wish to compute the
topological invariants of the Fourier-Mukai transform $\Cal
S_X(\Cal F)=R\hat\pi_\ast(\pi^\ast \Cal F{\overset L\to\otimes}\Cal Q)\in
D(\whatx)$ of a
sheaf $\Cal F$ on $X$ in
terms of those of $\Cal F$. The formula is obtained as usual
by the Riemann-Roch theorem, taking into account that we can compute the
Chern character
$\gamma$ of $\Cal Q$ from the sequence (5).

Let us define the Mukai vector and the Euler characteristic of the Fourier-Mukai
transform
$\Cal S_X(\Cal F)\in D(\whatx)$ by
$v(\Cal S_X(\Cal F))=\sum_{i=0}^2(-1)^i v(R^i\hat\pi_\ast(\pi^\ast\Cal
F\otimes\Cal Q))$ and $\chi(\Cal S_X(\Cal
F))=\sum_{i=0}^2(-1)^{i+1}\chi(R^i\hat\pi_\ast(\pi^\ast\Cal F\otimes\Cal Q))$.
 \proclaim{Proposition 11} Given a coherent sheaf $\Cal F$ on $X$,
let  $u=(\rho,c_1,\sigma)=(\rk(\Cal
F),\allowmathbreak c_1(\Cal F),\rk{\Cal F}+\ch_2(\Cal F))$ be the Mukai
vector of
$\Cal F$, and
$d=c_1\cdot H$. If $\hat u=v(\Cal S_X(\Cal F))=(\hat\rho,\hat c_1,\hat\sigma)$,
one has
$$
\align \hat\rho&=-3\rho+2\sigma+\ell\cdot c_1,\\
\hat c_1 &= (\ell\cdot
c_1+2d)\what H +(\rho+d-s)\hat\ell-\Psi^\ast(c_1), \\
\hat\sigma&=2\rho-3\sigma-\ell\cdot c_1\,.\endalign
$$
Then
$\chi(\Cal
S_X(\Cal F))=-\chi({\Cal F})$ and $\hat u^2=u^2$.
\qed\endproclaim
\proclaim{Corollary 8} The Fourier-Mukai
transform preserves the Euler characteristic and the degree of
WIT$_1$-sheaves, that is, if
$\Cal F$ is WIT$_1$ then $\chi(\hathsec F)=\chi(\Cal F)$ and $c_1(\Cal
F)\cdot H=
c_1(\what{\Cal F})\cdot\what H$.
\qed\endproclaim

\medskip\noindent{\bf Final remarks.}
Theorem 3 can be exploited to investigate the structure of the
entire moduli space of stable sheaves on a reflexive K3 surface $X$.
Let $u$ be a Mukai vector, $u\neq(2,-\ell,-3)$.
The locally free $\mu$-stable sheaves $\Cal F$ on $X$ with
$v(\Cal F)=u$ are a Zariski open set $M_H^\mu(u)\subset M_H(u)$; they are mapped
by the  Fourier-Mukai transform onto the Zariski open set of
locally free $\mu$-stable sheaves in $M_{\what H}(\hat u)$, where
$\hat u$ is given in terms of $u$ according to Proposition 11. This map
preserves the holomorphic symplectic structures of these spaces. 
Moreover, one has  $\dim M_H(u)=\dim M_{\what H}(\hat u)$ (since  $u^2=\hat
u^2$), and  so --- provided $M_H^\mu(u)$ is not empty --- there is  a
birational correspondence $M_H(u)\to M_{\what H}(\hat u)$.

In some cases   stronger results can be obtained; for instance it
can be shown that for any $n\ge 1$ the moduli space $M_{\what
H}(1+2n,-n\hat\ell,1-3n)$ is biholomorphic to the punctual Hilbert
scheme $\operatorname{Hilb}^n(X)$
\cite{8}.

In \cite{4} we give a completely algebraic proof of Theorem  3. 
Also in
\cite{4} we prove algebraically that $\O_{\whatx}(2\what H)$ is
the determinant line bundle, which is an
alternative proof of the ampleness of $\what H$.

The transcendental proof of Theorem 1 extends directly to higher dimensions,
providing a proof of the the fact that the Fourier-Mukai transform on
hyperk\"ahler manifolds maps quaternionic instantons to quaternionic 
instantons \cite{6}.

\medskip\noindent{\bf Acknowledgments.} We thank P\. Francia,
J\. M\. Mu\~noz Porras, K\. O'Grady and especially A\. Maciocia
for useful discussions and suggestions.
The first author also thanks S\. Donaldson for advice at a
preliminary stage of this work. This work was partly done while the
first author
was visiting the State University of New York at Stony Brook,
and the second
author was visiting Victoria University of Wellington, New Zealand.
They thank
the respective Departments of Mathematics for their warm hospitality
and for
providing excellent working conditions.

\obligedskip\obligedskip

\enddocument
\bye